\documentclass[11pt]{article} 
\title{{\bf On Orbit Equivalence of Quasiconformal Anosov Flows}}
\author{Yong Fang} 
\date{\it U.M.R. 
7501 du C.N.R.S, 
Institut de Recherche Math\'ematique Avanc\'ee, 7 rue 
Ren\'e Descartes, 67084 Strasbourg C\'edex, France\\ (e-mail : fang@math.u-strasbg.fr)} 
 \usepackage{amsmath,amsthm,amsfonts,amssymb,amscd,amstext} 
\chardef\bslash=`\\

\theoremstyle{definition}

\theoremstyle{remark}

\begin{document} 
\maketitle 
\renewcommand{\sectionmark}[1]{} 
{\bf Abstract} -- {\it We classify up to $C^\infty$ 
orbit equivalence the volume-preserving quasiconformal 
Anosov flows whose strong stable and strong 
unstable distributions are at least 
three-dimensional. If one of the strong distributions is two-dimensional, then we get a 
partial classification. Using these classification results, we obtain the following rigidity result :

Let $\Phi$ the the orbit foliation of the geodesic flow of a closed hyperbolic manifold of dimension at least three. 
Let $\Psi$ be another $C^\infty$ one-dimensional foliation. If $\Phi$ is $C^1$ conjuguate to $\Psi$, then 
$\Phi$ is $C^\infty$ conjuguate to $\Psi$.  
%of the 
%orbit foliations of the geodesic flows of 
%closed hyperbolic manifolds of dimension at least three.
}

$\ $

{\bf R\'esum\'e} -- {\it Nous classifions, \`a \'equivalence orbitale $C^\infty$ pr\`es, les flots d'Anosov quasiconformes 
topologiquement transitifs dont les distributions stable forte et instable forte sont de dimension au moins $3$. Si ces 
deux distributions sont de dimension au moins $2$, alors nous obtenons une classification partielle. Nous d\'eduisons 
de ces r\'esultats de classification le r\'esultat de rigidit\'e suivant :

Soit $\Phi$ le feuilletage orbitale du 
flot g\'eod\'esique d'une vari\'et\'e hyperbolique ferm\'ee de dimension au moins $3$. 
Soit $\Psi$ un autre feuilletage $C^\infty$ de dimension $1$. Si $\Phi$ et $\Psi$ sont $C^1$ conjugu\'e, alors 
$\Phi$ et $\Psi$ sont $C^\infty$ conjugu\'e.}\\\\ 
{\bf 1. Introduction}\\
{\bf 1.1. Motivation}

Let $M$ be a $C^{\infty}$-closed manifold. A $C^{\infty}$-flow $\phi_{t}$ generated 
by the non-singular vector field $X$ is said to be an {\it Anosov flow} if there exists a $\phi_t$-invariant 
splitting of the tangent bundle $$ TM= \mathbb{R}X\oplus E^{+}\oplus E^{-},$$
a Riemannian metric on $M$ and two positive numbers $a$ and $b$ such that 
for any $u^{\pm}\in E^{\pm}$ and for any $t \geq 0$, 
$$\parallel D\phi_{\mp t}(u^{\pm})\parallel\leq
a e^{-bt}\parallel u^{\pm}\parallel ,$$
where $E^{-}$ and $E^{+}$ are said to be the {\it strong stable} and {\it strong unstable} distributions of the flow. 
For any $x\in M$ the leaves containing $x$ of $E^+$ and $E^-$ are denoted respectively by 
$W^+_x$ and $W^-_x$.

Define two functions on $M\times\mathbb{R}$ as following:
$$K^+(x, t) = \frac{ max \{ \parallel D\phi_t(u)\parallel \  \mid u\in E^+_x,\  \parallel u\parallel =1\}}
{min \{ \parallel D\phi_t(u)\parallel \  \mid u\in E^+_x, \  \parallel u\parallel =1\}}$$
and 
$$K^-(x, t) = \frac{ max \{ \parallel D\phi_t(u)\parallel \  \mid u\in E^-_x,\  \parallel u\parallel =1\}}
{min \{ \parallel D\phi_t(u)\parallel \  \mid u\in E^-_x, \  \parallel u\parallel =1\}}.$$
If $K^-$ $(K^+)$ is bounded, then the Anosov flow $\phi_t$ is said to be 
{\it quasiconformal on the stable (unstable) distribution}. If $K^+$ and $K^-$ 
are both bounded, then $\phi_t$ is said to 
be {\it quasiconformal}. If it is the case, then the superior bound of $K^+$ and $K^-$ is 
said to be the {\it distortion} of $\phi_t$. 
The corresponding notions for Anosov diffeomorphisms are defined similarly (see {\bf [Sa]}).

Recall that two $C^\infty$ Anosov flows $\phi_t: M\to M$ and $\psi_t: N\to N$ are said to be 
{\it $C^k$ flow equivalent} $(k\geq 0)$ if 
there exists a $C^k$ diffeomorphism $h: M\to N$ such that $\phi_t= h^{-1}\circ \psi_t\circ h$ for all $t\in\mathbb{R}$. 
They are said to be {\it $C^k$ orbit equivalent} $(k\geq 0)$ if there exists a 
$C^k$ diffeomorphism $h: M\to N$ sending the orbits of $\phi_t$ onto the orbits of $\psi_t$ such that 
the orientations of the orbits are preserved. Similarly, two $C^\infty$ foliations 
$\Phi$ and $\Psi$ are said to be {\it $C^k$ conjuguate} if there exists a $C^k$ diffeomorphism 
$h: M\to N$ sending the leaves of $\Phi$ onto those of $\Psi$. 
%the flow $\chi_t= h^{-1}\circ \psi_t\circ h$ is a time change of $\phi_t$. 
By convention, a $C^0$ diffeomorphism means a homeomorphism. We can prove easily the following\\\\
{\bf Lemma 1.1.} {\it Let $\phi_t$ and $\psi_t$ be two $C^\infty$ Anosov flows. 
If they are $C^1$ orbit equivalent and 
$\psi_t$ is quasiconformal, then $\phi_t$ is also quasiconformal.}\\\\
{\bf Proof.} Denote by $\phi: M\to N$ the $C^1$ orbit conjugacy 
between $\phi_t$ and $\psi_t$ and by $\bar{\mathcal F}^\pm$ and 
$\bar{\mathcal F}^{\pm, 0}$ the Anosov foliations of $\psi_t$. Then we have the 
following\\\\
{\bf Sublemma.} {\it Under the notations above, $\phi({\mathcal F}^{\pm, 0})= \bar{\mathcal F}^{\pm, 0}$.}\\\\
{\bf Proof.} Define $\hat\phi_t= \phi\circ\phi_t\circ \phi^{-1}$. Then $\hat\phi_t$ is a $C^1$ flow on $N$ with the same 
orbits as $\psi_t$. So there exists a $C^1$ map $\alpha: \mathbb R\times N\to \mathbb R$ such that 
$\hat\phi_t(x)= \psi_{\alpha(t, x)}(x)$. Define $\hat E^-= (D\phi)(E^-)$. Then it is the $C^0$ 
tangent bundle of the $C^1$ foliation $\hat{\mathcal F}^-= \phi({\mathcal F}^-)$. 

Let us prove first that $\hat E^-\subseteq \bar E^{-, 0}$. Fix a $C^0$ Riemannian metric $g$ on $N$ such 
that $\bar E^+$ and $\bar E^-$ and $\bar X$ are orthogonal to each other. Since $\phi$ is $C^1$, then it is 
bi-Lipschitz. We deduce that for all $x\in N$ and $\hat u\in \hat E^-_x$, 
$\parallel D\hat\phi_t(\hat u)\parallel\to 0$ if $t\to +\infty$. 

If $\hat u= \bar u^++ a\bar X_x+ \bar u^-$ and $\bar 
u^+\not =0$, then by a simple calculation we get for a certain function 
$b_x: \mathbb R\to \mathbb R,$
$$ D\hat\phi_t(\hat u)= D(\psi_{\alpha(t, x)})(\bar u^+)+ b_x(t)\cdot \bar X_{\hat\phi_t(x)}+D(\psi_{\alpha(t, x)})(\bar u^-).$$
So we get $\parallel D\hat\phi_t(\hat u)\parallel\geq  \parallel D(\psi_{\alpha(t, x)})(\bar u^+)\parallel\to +\infty$ if 
$t\to +\infty$, which is a contradiction. We deduce that $\hat E^-\subseteq \bar E^{-, 0}$. So $\phi({\mathcal F}^0)
\subseteq\bar {\mathcal F}^{-, 0}$. Then it is easy to see that $\phi({\mathcal F}^{-,0})=
\bar{\mathcal F}^{-, 0}$, i.e. $\phi$ sends $C^1$ diffeomorphically 
each leaf of ${\mathcal F}^{-,0}$ onto a leaf of $\bar{\mathcal F}^{-, 0}$. Similarly we have 
$\phi({\mathcal F}^{+,0})=
\bar{\mathcal F}^{+, 0}$.  $\square$

$\ $

In particular we deduce from the sublemma above that 
$ {\hat E}^\pm\oplus\mathbb R\bar X= {\bar E}^\pm\oplus \mathbb R\bar X$. 
By projecting parallel to the direction of $\bar X$, we get a 
$C^0$ section $P$ of $End({\hat E}^+,{\bar E}^+)$ and two positive constants $A_1$ and $A_2$ such that 
$$ A_1\parallel u\parallel \leq \parallel P(u)\parallel\leq A_2\parallel u\parallel,\  \forall\  u\in {\hat E}^+.$$
Then it is easy to verify that for all $x\in N$, $P\circ D_x\hat\phi_t= D_x\psi_{\alpha(t, x)}\circ P$.  

If $\psi_t$ is quasiconformal with 
distortion $K$, then we have the following estimation
$$\hat K^+(x, t)=\frac{sup\{ \parallel D\hat\phi_t(\hat u^+)
\parallel \  \mid\   \parallel \hat u^+\parallel =1,\  \hat u^+\in\hat E^+_x\}}
{inf \{ \parallel D\hat\phi_t(\hat u^+)\parallel \  \mid \  \parallel \hat u^+\parallel =1,\  \hat u^+\in\hat E^+_x\}}$$
$$\leq 
\frac{A_2}{A_1}\cdot \frac{sup\{
\parallel P(D\hat\phi_t(\hat u^+))\parallel \  \mid\   \parallel \hat u^+\parallel =1,\  \hat u^+\in\hat E^+_x\}}
{inf \{ \parallel P(D\hat\phi_t(\hat u^+))\parallel \  \mid\   \parallel \hat u^+\parallel =1,\  \hat u^+\in\hat E^+_x\}}$$
$$\leq (\frac{A_2}{A_1})^2\cdot \frac{sup\{\parallel D_x\psi_{\alpha(t, x)}(\frac{
P(\hat u^+)}{\parallel P(\hat u^+)\parallel})\parallel \  \mid\   \parallel \hat u^+\parallel =1,\  \hat u^+\in\hat E^+_x\}}
{inf \{ \parallel D_x\psi_{\alpha(t, x)}(\frac{
P(\hat u^+)}{\parallel P(\hat u^+)\parallel}) \parallel \  
\mid\   \parallel \hat u^+\parallel =1,\  \hat u^+\in\hat E^+_x\}}$$
$$\leq (\frac{A_2}{A_1})^2\cdot K.$$
Similarly $\hat K^-$ is also bounded. So $\hat\phi_t$ is quasiconformal. 
Since $\phi$ is a $C^1$ diffeomorphism, then it is bi-Lipschitz. We deduce that $\phi_t$ is also quasiconformal. 
$\square$    

$\ $

It is easy to see that the geodesic flow of a closed hyperbolic manifold is quasiconformal (even 
conformal). Then by the previous lemma, each $C^\infty$ time change of such a flow is also quasiconformal. 
It is easily seen that an Anosov diffeomorphism is quasiconformal 
iff its suspension is quasiconformal. So if $\phi$ 
denotes a semisimple hyperbolic automorphism of a torus with two enginvalues, then its 
suspension is a quasiconformal Anosov flow. 

It seems to be a common phenomena in mathematics that things can only be effectively studied and 
understood when placed in a suitable and flexible environment. Conformal structures (Anosov flows) are 
pretty rigid while quasiconformal structures (Anosov flows) seem to be much more 
flexible. We wish to better understand the classical conformal 
Anosov flows, notably the geodesic flows of closed hyperbolic manifolds, by using quasiconformal 
techniques, which is our motivation to study general quasiconformal Anosov systems.    \\\\
{\bf 1.2. Main theorems}

In our previous paper {\bf [Fa1]}, we have studied the rigidity of volume-preserving 
Anosov flows with smooth $E^+\oplus E^-$. In particular we have obtained the following\\\\
{\bf Theorem 1.1.} ({\bf [Fa1]}, Corollary $1$) {\it Let $\phi$ be a $C^\infty$ volume-preserving quasiconformal Anosov 
diffeomorphism on a closed manifold 
$\Sigma$. If the dimensions 
of $E^+$ and $E^-$ are at least two, then up to finite covers 
$\phi$ is $C^\infty$ conjugate to a hyperbolic automorphism of a torus.}

$\ $

In {\bf [K-Sa]}, B. Kalinin and V. Sadovskaya 
have classified the topologically transitive quasiconformal Anosov diffeomorphisms whose strong stable and 
unstable distributions are of dimension at least $3$. Their argument, though quite elegant, meets 
an essential difficulty in the case that one of the strong stable and unstable distributions is 
two dimensional.  

Based on the classification result 
in {\bf [K-Sa]}, we classify in this paper completely the quasiconformal Anosov flows whose strong 
stable and unstable distributions have relatively high dimensions. More precisely, we prove\\\\ 
{\bf Theorem 1.2.} {\it Let $\phi_t$ be a $C^\infty$ 
topologically transitive quasiconformal Anosov flow such that 
$E^+$ and $E^-$ are at least three dimensional. Then 
up to finite covers, $\phi_t$ is $C^\infty$ orbit equivalent either to the geodesic flow of a hyperbolic 
manifold or to the suspension of a hyperbolic automorphism of a torus.

Under the conditions above, if $E^+\oplus E^-$ is in addition $C^1$, then up to a constant 
change of time scale and finite covers, $\phi_t$ is $C^\infty$ flow equivalent either to a 
canonical time change of the geodesic flow of a hyperbolic manifold or to the 
suspension of a hyperbolic automorphism of a torus.}

$\ $

Recall that if $\alpha$ is a closed $C^\infty$ $1$-form on $M$ such that $1+\alpha(X)>0$, then 
the flow of $\frac{X}{1+\alpha(X)}$ is said to be a {\it canonical time change} of $\phi_t$, 
where $X$ denotes the generator of $\phi_t$. 

If one of the strong distributions is two-dimensional, then by using our Theorem 
$1.1$, we get the following partial result.\\\\ 
{\bf Theorem 1.3.} {\it Let $\phi_t$ be a $C^\infty$ 
volume-preserving quasiconformal Anosov flow such that 
$E^+$ is of 
dimension two and $E^-$ is of dimension at least two. 
If $\phi_t$ has the sphere-extension property, then 
up to finite covers, $\phi_t$ is $C^\infty$ orbit equivalent either to the geodesic flow of a 
three-dimensional hyperbolic manifold or to the suspension of a hyperbolic automorphism of a torus.

Under the conditions above, if $E^+\oplus E^-$ is in addition $C^1$, then up to a constant change of time scale 
and finite covers, $\phi_t$ is $C^\infty$ flow equivalent either to a canonical time change 
of the geodesic flow of a three-dimensional hyperbolic manifold or to the suspension of a 
hyperbolic automorphism of a torus. }

$\ $

The sphere-extension property will be defined in Subsection $2.3$. Let us just mention that this property is 
invariant under $C^1$ orbit equivalence.

Now we can get some concrete applicaions of our results above. Recall at first 
that flow conjugacies between Anosov flows have been and being extensively studied. 
The philosophical conclusion is that they exist rarely, even $C^0$ ones. Let us just mention two of the most 
beautiful supporting results (see also {\bf [L1]} and {\bf [L2]}) :\\\\
{\bf Theorem 1.4.} (U. Hamenst$\ddot{\rm a}$dt, {\bf [Ham2]}) {\it Let $M$ be 
a closed negatively curved manifold. If the geodesic flow of $M$ is $C^0$ 
flow equivalent to that of a locally symmetric space of rank one $N$, then $M$ is isometric to $N$.}\\\\
{\bf Theorem 1.5.} (R. de la Llave and R. Moriyon, {\bf [LM]}) {\it Let $\phi_t$ and $\psi_t$ 
be two $C^\infty$ three-dimensional 
volume-preserving Anosov flows. If they are $C^1$ flow equivalent, then they are 
$C^\infty$ flow equivalent.}

$\ $

However there exist plenty of $C^0$ orbit 
conjugacies between Anosov flows. For example, if two $C^\infty$ 
Anosov flows are sufficiently $C^1$-near, then they are H$\ddot{\rm o}$lder-continuous orbit equivalent by the 
celebrated structural stability (see {\bf[An]}). A natural 
question to ask is whether 
$C^1$ orbit conjugacies between Anosov flows are rare. 

We can deduce from Theorems $1.1$ and $1.2$ the 
following result showing that $C^1$ orbit conjugacies are surely rare 
in some cases, while H$\ddot{\rm o}$lder-continuous orbit conjugacies 
are abundant.\\\\
{\bf Theorem 1.6.} {\it Let $\phi_t$ be a $C^\infty$ Anosov flow and $\psi_t$ be the geodesic flow of a closed 
hyperbolic manifold of dimension at least three. If 
$\phi_t$ and $\psi_t$ are $C^1$ orbit equivalent, then they are $C^\infty$ orbit equivalent.}

$\ $

In order to state the next result, let us recall firstly some notions.\\\\
{\bf Definition 1.1.} Let $(X, d_X)$ and $(Y, d_Y)$ be two metric spaces. Then they are 
said to be {\it quasi-isometric} if there a map $f: X\to Y$ and two positive numbers $C$ and $D$ such that the 
following two conditions are satisfied:\\
$(1)$ $\frac{1}{C} d_X(x, y)-D\leq d_Y(f(x), f(y))\leq C d_X(x, y)+D, \  \forall \  x, y\in X.$\\
$(2)$ For any $y\in Y$, there exists $x\in X$ such that $d_Y(y, f(x))\leq D$.

$\ $

Roughly speaking, two metric spaces are quasi-isometric if and only if they are 
bi-Lipschitz equivalent in the large scale. 
Recall that for any $n\geq 2$, a $n$-dimensional Riemannian manifold $M$ is said to be {\it hyperbolic} if it 
has constant sectional curvature $-1$. We 
denote by $\mathbb H^n$ the unique simply connected hyperbolic manifold. 
By combining some classical results 
with our previous theorem, we get the following\\\\
{\bf Corollary 1.1.} {\it Let $M$ be a $n$-dimensional closed Riemannian manifold of negative curvature 
such that $n\geq 3$. Then we have the following relations between dynamics and geometry:\\
$(1)$ The geodesic flow of $M$ is H$\ddot{ o}$lder-continuously orbit equivalent to that of a hyperbolic manifold 
if and only if the universal covering 
space $\widetilde M$ with its lifted metric is quasi-isometric to $\mathbb H^n$.\\
$(2)$ The geodesic flow of $M$ is $C^1$ orbit equivalent to that of a hyperbolic manifold 
if and only if $M$ has constant negative curvature.}

$\ $

Let $g$ be the Riemannian metric of a closed hyperbolic manifold of dimension at least three. 
Let $g'$ be a perturbed Riemannian metric of $g$ with non-constant negative curvature. Thus by the proof of the 
previous corollary, the orbit foliation of the geodesic flow of $g$ is $C^0$ conjuguate to that of $g'$. However, 
by the previous corollary, these two $C^\infty$ one-dimensional foliations are not 
$C^1$ conjuguate.

By combining a result of R. Man$\tilde{\rm e}$ with our Theorem $1.6$, we get finally 
the following corollary, which is 
the key objective of this article.\\\\
{\bf Corollary 1.2.} {\it Let $\Phi$ be the orbit foliation of the geodesic flow of a 
closed hyperbolic manifold of dimension at least 
three. Let $\Psi$ be another $C^\infty$ one-dimensional foliation. If $\Phi$ and $\Psi$ are $C^1$ conjuguate, then 
they are $C^\infty$ conjuguate.}\\\\
{\bf 1.3. The organization of the paper}

In Section 
two we recall some fundamental facts concerning quasiconformal Anosov flows 
and some properties of 
transverse $(G, T)$-structures of foliations. Then in Section three we 
prove Theorems $1.2$ and $1.3$. Finally in Section four we apply these two results 
to geodesic flows to deduce Theorem $1.6$ and Corollaries $1.1$ and $1.2$.\\\\  
{\bf 2. Preliminaries}\\
{\bf 2.1. Linearizations and smooth conformal structures}

Let us recall firstly the following two results established in {\bf [Sa]}:\\\\
{\bf Theorem 2.1.} ({\bf [Sa]}, Theorem $1.3$) {\it Let $f$ be a topologically 
transitive $C^\infty$ Anosov diffeomorphism ($\phi_t$ be a topologically mixing $C^\infty$ Anosov flow) on a closed 
manifold $M$ which is quasiconformal on the unstable distribution. Then it is conformal with respect to 
a Riemannian metric on this distribution which is continuous on $M$ and $C^\infty$ along the leaves of the unstable 
foliation.}\\\\
{\bf Theorem 2.2.} ({\bf [Sa]}, Theorem $1.4$) {\it Let $f$ ($\phi_t$) be a $C^\infty$ Anosov 
diffeomorphism (flow) on a closed manifold $M$ with dim$E^+\geq 2$. Suppose that it is conformal with respect to 
a Riemannian metric on the unstable distribution which is continuous on $M$ and $C^\infty$ along the leaves of 
the unstable foliation. Then the (weak) stable holonomy maps are conformal and the (weak) stable distribution is 
$C^\infty$.}

$\ $

Let us recall briefly the steps to prove these two theorems in the case of flow. Denote by $\phi_t$ a 
topologically mixing quasiconformal Anosov flow. Then by some classical arguments (see 
{\bf [Su]} and {\bf [Tu]}), V. Sadovskaya 
found two measurable $\phi_t$-invariant conformal structures $\bar\tau^+$ and 
$\bar\tau^-$ along respectively the leaves of ${\mathcal F}^+$ and 
${\mathcal F}^-$. Then as is usual for Anosov flows, these two 
conformal structures were pertubated to continuous $\phi_t$-invariant ones denoted by $\tau^+$ and $\tau^-$. 
Using the linearizations, she proved that along each leaf of ${\mathcal F}^+$, $\tau^+$ is isometric 
to a vector space with its canonical conformal struture, which has permitted her to 
blow up the smoothness of weak holonomy maps. Then by using a result of J. L. Journ\'e, she proved 
the smoothness of the weak stable and unstable distributions.

In {\bf [Fa1]}, we have 
proved the following lemma based on {\bf [Pl1]}. For the sake of completeness, let us recall the arguments.\\\\
{\bf Lemma 2.1.} {\it Let $\phi_t$ be a $C^\infty$ topologically transitive Anosov flow. Then 
we have the following alternative:\\
$(1)$ $\phi_t$ is topologically mixing,\\ 
$(2)$ $\phi_t$ admits a $C^\infty$ closed global section with constant return time.}\\\\
{\bf Proof.} If Case $(2)$ is true, then up to a constant change of time scale, $\phi_t$ is $C^\infty$ flow equivalent to 
the suspension of a 
$C^\infty$ Anosov diffeomorphism. Thus it is not topologically mixing. So the alternative is exclusive.

If there exists $x\in M$ such that ${W^+_x}$ is not dense in $M$, then 
by Theorem $1.8$ of {\bf [Pl1]} $E^+\oplus E^-$ is the tangent bundle 
of a $C^1$ foliation $\mathcal F$. In addition the leaves of $\mathcal F$ are all comact. So $\phi_t$ admits a $C^1$ 
closed global section with constant return time. Then to realize Case $(2)$ we need only prove that 
$E^+\oplus E^-$ is $C^\infty$.
  
Denote by $\lambda$ the canonical $1$-form of $\phi_t$. Then $\lambda$ is, a priori, a continuous $1$-form on $M$. 
For each point $y\in M$ we take a small neighborhood $F_y$ of $y$ in the leaf 
containing $y$ of $\mathcal F$. Then we can construct a local $C^1$ chart 
$\theta_y: (-\epsilon, \epsilon)\times F_y\to M$ such that $
\theta_y(t, z)= \phi_t(z)$. In this chart we have $\lambda= dt.$ We deduce that $\int_{\gamma} 
\lambda=0$ for each piecewise $C^1$ closed curve $\gamma$ contained in the image of $\theta_y$. So $\lambda$ 
is locally closed (see Section two of {\bf [Pl1]} for the definition). 
Then by Proposition $2.1$ of {\bf [Pl1]}, $\lambda$ is seen to be closed in a 
weak sense, i.e. for every $C^1$ immersed two-disk 
$\sigma$ such that $\partial\sigma$ is piecewise $C^1$, 
$$\int_{\partial\sigma}\lambda=0.$$
So by integrating along the closed curves, $\lambda$ gives an element in $Hom(\pi_1(M), \mathbb R)$, i.e. the space 
of group homomorphisms of $\pi_1(M)$ into $\mathbb R$, 
where $\pi_1(M)$ denotes the 
fundamental group of $M$. However we have naturally
$$Hom(\pi_1(M), \mathbb R)\cong (H_1(M, \mathbb R))^\ast\cong H^1(M, \mathbb R),$$
where $H^1(M, \mathbb R)$ denotes the first de Rham cohomology group of $M$. 
So there exists a $C^\infty$ $1$-form $\beta$ such that for each $C^\infty$ closed curve $\gamma$, 
$$\int_\gamma \lambda=\int_\gamma \beta.$$
So by integrating $(\lambda-\beta)$ along curves, we get on $M$ a well-defined 
continuous function $f$. Then for any $y\in M$ and any $t\in\mathbb R$ we have
$$ (f\circ\phi_t)(y) - f(y) = \int_0^t(1-\beta(X))(\phi_s(y)) ds.$$

Since the right-hand side of this identity is a $C^\infty$ $\mathbb{R}$-cocycle, then by {\bf [LMM]} 
$f$ is smooth. However by the definition of $f$ we have 
$$\lambda-\beta= df.$$ 
So $\lambda$ is also smooth. We deduce that $E^+\oplus E^-(= Ker\lambda)$ is $C^\infty$. 
So Case $(2)$ is realized if $W^+_x$ is not dense 
for a certain point $x\in M$.

Now suppose that  for all $x\in M$, $W^+_x$ is dense in $M$. Fix a Riemannian metric on $M$. For 
any $x\in M$ and any $r >0$ we denote by $M_{x,r}$ and $W_{x,r}^{+}$ the balls of 
center $x$ and radius $r$ in $M$ and 
$W^+_x.$ Take arbitrarily two open subsets $U$ and $V$ in $M$ and a small ball $M_{y,\epsilon}$ in $V$. Since 
$M$ is closed and each strong unstable leaf is supposed to be dense in $M$, then we can find $R < +\infty$ such that 
$$W^{+}_{x,R} \cap M_{y,\epsilon}\not =\emptyset, \  \forall \  x\in M.$$
Take a small disk $W_{x,\delta}^{+}$ in $U$. Then by the Anosov property there exists $T >0$ such that 
$$\phi_t( W_{x,\delta}^{+})\supseteq W_{\phi_t(x),R}^{+}, \  \forall \  t\geq T.$$   
So $\phi_tU\cap V\not =\emptyset, \  \forall \  t\geq T,$ i.e. $\phi_t$ is 
topologically mixing.  $\square$

$\ $

If $\phi_t$ is a topologically transitive quasiconformal Anosov flow 
such that the dimensions of $E^+$ and $E^-$ are at least two, then by the previous lemma 
and Theorem $2.1$, it preserves two continuous conformal strutures $\tau^+$ and $\tau^-$ which are $C^\infty$ 
along the leaves of ${\mathcal F}^+$ and ${\mathcal F}^-$. Then by Theorem $2.2$, these two 
conformal structures are invariant under the weak holonomy maps. 
In addition we can deduce from the previous lemma 
and Theorem $2.2$ the following\\\\  
{\bf Lemma 2.2.} {\it Let $\phi_t$ be a $C^\infty$ topologically transitive quasiconformal Anosov flow 
such that the dimensions of $E^+$ and $E^-$ are at least two. Then $E^{+, 0}$ and $E^{-, 0}$ are both $C^\infty$.}

$\ $

Based on {\bf [Sa]}, we have found in {\bf [Fa1]} for each quasiconformal Anosov flow $\phi_t$ an 
{\it unstable linearization} $\{ h^+_x\}_{x\in M}$ such that for any $x\in M$, 
$h^+_x: W^+_x\to E^+_x$ is a $C^\infty$ diffeomorphism and the following 
conditions are satisfied: \\
$(1)$ $h^{+}_{\phi_s(x)}\circ \phi_s = D_x\phi_s\circ h^{+}_x,\  \forall\  s\in\mathbb R$,\\
$(2)$ $h^{+}_x(x) =0$ and $(Dh^{+}_x)_x$ is the identity map,\\
$(3)$ $h^{+}_x$ depends continuously on $x$ in the $C^\infty$ topology.

Recall that the family of diffeomorphisms satisfying these three conditions 
is unique. Similarly we have the {\it stable linearization } $\{ h^-_x\}_{x\in M}$ of $\phi_t$.

Suppose that $\phi_t$ 
satisfies the conditions of Lemma $2.2$. For any $ x\in M$ we can extend the conformal structure $\tau^+_x$ 
at $0\in E^+_x$ to all other points of $E^+_x$ via linear translations. If 
the resulting translation-invariant conformal 
structure on $E^+_x$ is denoted by $\sigma^+_x$, then $(E^+_x, \sigma^+_x)$ is isometric to the 
canonical conformal structure of $\mathbb R^n$ if $E^+$ is $n$-dimensional. Similarly we have 
$(E^-_x, \sigma^-_x)$ for any $x\in M$.
  
Then by Lemma $3.1$ of {\bf [Sa]} we get for any $x\in M$, 
$(h^+_x)_\ast \tau^+\mid_{W^+_x}=\sigma^+_x$ and $(h^-_x)_\ast \tau^-\mid_{W^-_x}=\sigma^-_x.$\\\\
{\bf 2.2. Transverse (G, T)-structures}

In this subsection, we consider the transverse $(G, T)$-structures of foliations. 
Let $\mathcal{F}$ be a $C^\infty$ foliation on a connected manifold $M$. Denote 
by $Q_{\mathcal F}$ the 
leaf space of $\mathcal F$ and by $\mathcal{F}(A)$ the saturation of 
$\mathcal{F}$ on $A$ for any $A\subseteq M$. We assume that the holonomy maps of $\mathcal{F}$ are 
defined on connected transverse sections.

Let $G$ be a real Lie group acting effectively and transitively on a connected manifold $T$. If $\Sigma$ is a $C^\infty$ 
transverse section of $\mathcal F$ and $\phi$ is a $C^\infty$ diffeomorphism of $\Sigma$ onto its open image in 
$T$, then $(\Sigma, \phi)$ is said to be a {\it transverse $T$-chart}. Two transverse 
$T$-charts $(\Sigma_1,\phi_1)$ and 
$(\Sigma_2,\phi_2)$ are said to be {\it compatible} if for each holonomy map 
$h$ of a 
germ of $\Sigma_1$ to a germ of $\Sigma_2$, the map $\phi_1\circ h\circ \phi_2^{-1}$ is 
locally the restriction of elements of $G$.

A family of transverse sections is said to be {\it covering} if each leaf of $\mathcal F$ intersects at least one of 
the sections in this family. By definition, a {\it transverse $(G, T)$-structure} on $\mathcal F$ is a maximal 
family of compatible transverse $T$-charts of which the underlying family of transverse sections is covering. 

In order to define a transverse $(G, T)$-structure, we need 
just seperate out a family of covering 
compatible $T$-charts. Then by considering all the $T$-charts compatible with this family, we get 
automatically a 
transverse $(G, T)$-structure.   

Denote by $\widetilde {\mathcal F}$ the lifted foliation on the 
universal covering space $\widetilde M$ of $M$ and 
denote by $\pi$ the projection of $\widetilde M$ onto $M$. For each transverse $(G, T)$-structure on $\mathcal F$, 
we get naturally a lifted transverse $(G, T)$-structure on $\widetilde{\mathcal F}$ 
by considering the composition of $\pi$ with the $T$-charts 
of the given transverse $(G, T)$-structure on $\mathcal F$. Then by {\bf [Go]} there exists a $C^\infty$ 
submersion $\mathcal D: \widetilde M\to T$ and a group homomorphism $\mathcal H: \pi_1(M)\to 
G$ satisfying the following two conditions:\\
$(1)$ ${\mathcal D}(\gamma x)= {\mathcal H}(\gamma){\mathcal D}(x),$ $\forall$ $x\in\widetilde M$, 
$\forall$ $\gamma\in \pi_1(M).$\\
$(2)$ The lifted foliation $\widetilde {\mathcal F}$ is defined by $\mathcal D$.\\
This submersion $\mathcal D$ is said to be the {\it developing map} of the transverse $(G, T)$-structure of 
$\mathcal F$ and $\mathcal H$ is said to be the {\it holonomy representation} of $\mathcal D$. 
The transverse $(G, T)$-structure of $\mathcal F$ is said to be {\it complete} if 
$\mathcal D$ is a $C^\infty$ fibre bundle over $\mathcal D(\widetilde M)$. 

If $\mathcal D'$ denotes another developing map with holonomy representation $\mathcal H'$, then by 
{\bf [Go]} there exists a unique element $g\in G$ such that 
$\mathcal D'= g\circ \mathcal D$ and $\mathcal H'= g\cdot \mathcal H\cdot g^{-1}$. 

Since $\mathcal D$ is obtained by analytic continuation along curves (see {\bf [Go]}), then for each transverse 
section $\Sigma$ of $\widetilde {\mathcal F}$ such that 
${\mathcal D}\mid_{\Sigma}$ is a $C^\infty$ diffeomorphism onto its image, ${\mathcal D}\mid_{\Sigma}$ is a 
transverse $T$-chart of the lifted transverse $(G, T)$-structure of $\widetilde{\mathcal F}$.

Since $\widetilde {\mathcal F}$ is defined by the submersion $\mathcal D$, then $\mathcal D$ sends 
each leaf of $\widetilde F$ to a point of $T$. 
Thus for any $x\in\widetilde M$ there exists a small $C^\infty$ transverse section $\Sigma$ containing $x$ of 
$\widetilde{\mathcal F}$ such that each leaf of $\widetilde {\mathcal F}$ intersects 
$\Sigma$ at most once. 
Then it is easily seen that each leaf of $\widetilde{\mathcal F}$ is closed in $\widetilde M$.  

Denote by $Q_{\widetilde {\mathcal F}}$ the leaf space of $\widetilde {\mathcal F}$. Then we have the quotient map 
${\bar {\mathcal D}}:Q_{\widetilde {\mathcal F}}\to \mathcal D(\widetilde M)$. Since each leaf of $\widetilde {\mathcal F}$ is closed, then 
$\bar {\mathcal D}$ is bijective iff the $\mathcal D$-inverse image of each point of $T$ is connected. 
If this is the case, then by 
considering the projections of small transverse $T$-charts, $Q_{\widetilde {\mathcal F}}$ 
becomes naturally a $C^\infty$ (seperable) manifold such that 
$\bar {\mathcal D}$ is a $C^\infty$ diffeomorphism of $Q_{\widetilde {\mathcal F}}$ onto ${\mathcal D}(\widetilde M)$. 
In addition, the 
fundamental group $\pi_1(M)$ of $M$ acts 
naturally on $Q_{\widetilde {\mathcal F}}$.

By {\bf [Hae]} we have the following \\\\
{\bf Proposition 2.1.} {\it Let $(M_1,\mathcal F_1)$ and $(M_2,\mathcal F_2)$ be two $C^\infty$ foliations 
with complete transverse $(G, T)$-structures. Suppose that their 
developing maps have both connected fibres and the holonomy covers of their leaves are all contractible. If the 
$\pi_1(M_1)$-action on $Q_{\widetilde {{\mathcal F}_1}}$ is $C^\infty$ conjugate to the 
$\pi_1(M_2)$-action on $ Q_{\widetilde{{\mathcal F}_2}}$, then there exists a $C^\infty$ map 
$f: M_1\to M_2$ such that the following conditions are satisfied:\\
$(1)$ $f$ is a surjective homotopy equivalence.\\
$(2)$ $f$ sends each leaf of ${\mathcal F}_1$ onto a leaf of ${\mathcal F}_2$ and $f$ sends different leaves to 
different leaves.\\
$(3)$ $f$ is transversally a local $C^\infty$ diffeomorphism conjugating the two transverse $(G, T)$-structures.}

$\ $

The lemma below is self-evident and will be used several times in the following.\\\\
{\bf Lemma 2.3.} {\it Let $(M, \mathcal F)$ be a $C^\infty$ foliation with 
a transverse $(G, T)$-structure. If $T_1$ is an open subset of $T$ and $G_1$ is a closed 
Lie subgroup of $G$ acting transitively on $T_1$ such that $\mathcal D(\widetilde M)\subseteq T_1$ and 
$\mathcal H(\pi_1(M))\subseteq G_1$, then $\mathcal F$ admits a transverse $(G_1, T_1)$-structure with the same 
developing map $\mathcal D$ and the same holonomy representation $\mathcal H$, 
which is compatible with the initial transverse $(G, T)$-structure.}

$\ $

By adapting the arguments in {\bf [Gh2]} we can prove the following \\\\
{\bf Lemma 2.4.} {\it Let $(M, \mathcal F)$ be a $C^\infty$ foliation with a transverse 
$(G, T)$-structure. If the closed leaves of $\mathcal F$ are dense in $M$ 
and the $\mathcal D$-inverse image of each point of $T$ is connected, then $\mathcal H(\pi_1(M))$ is a discrete 
subgroup of $G$.}\\\\
{\bf  Proof.} Denote by $\pi_1$ the projection of $\widetilde M$ onto $Q_{\widetilde{\mathcal F}}$ 
and by $\pi_2$ the projection of $Q_{\widetilde{\mathcal F}}$ onto  $Q_{\mathcal F}$ the 
leaf space of $\mathcal F$. Take a closed leaf 
$F_x$ of $\mathcal F$ and $\widetilde x\in\widetilde M$ such that $\pi(\widetilde x)=x$. Since 
$F_x$ is closed, then we can find a fine transverse 
section $\Sigma$ passing through $\widetilde x$ such that $\pi$ sends $\Sigma$ diffeomorphically onto its 
image and $ F_x\cap \pi(\Sigma)= \{x\}.$ 
So for each $y\in \Sigma$ and $y\not = \widetilde x$, $\pi(y)$ is not in $F_x$. Thus $\pi_2^{-1}(F_x)$ is discrete 
in $Q_{\widetilde{\mathcal F}}$.

Since $F_x$ is closed and $\pi^{-1}(F_x)= \pi_1^{-1}(\pi_2^{-1}(F_x)),$ 
then $\pi_2^{-1}(F_x)$ is closed in $\widetilde Q_{\mathcal F}$. 
So the $\pi_1(M)$-orbit of $\widetilde F_{\widetilde x}$, i.e. $\pi_2^{-1}(F_x)$ is closed and discrete in 
$Q_{\widetilde{\mathcal F}}$. 

Since the closed leaves of $\mathcal F$ are dense in $M$ i.e. the union of 
all the closed leaves is dense, then $\pi_2$-inverse 
images of these closed leaves form a dense subset $P$ of $Q_{\widetilde{\mathcal F}}$ such 
that the $\pi_1(M)$-orbit of each point of $P$ is closed and discrete.

Suppose on the contrary that $\mathcal H(\pi_1(M))$ is not discrete in $G$. Since the 
$\mathcal D$-inverse image of each point of $T$ is connected, then $\mathcal D$ induces a $C^\infty$ 
diffeomorphism 
$\bar{\mathcal D}:Q_{\widetilde{\mathcal F}}\to {\mathcal D}(\widetilde M)$. So $\bar{\mathcal D}(P)$ is dense in 
$\mathcal D(\widetilde M)$ and the $\mathcal H(\pi_1(M))$-orbit of each point of $\bar{\mathcal D}(P)$ 
is discrete and closed in $\mathcal D(\widetilde M)$. 

Take a non-trivial one-parameter subgroup $g_t$ of the closure of $\mathcal H(\pi_1(M))$ in $G$. For each $t
\in \mathbb R$, $g_t$ preserves the closed complement of $\mathcal D(\widetilde M)$. So we have 
$$g_t(\mathcal D(\widetilde M))= \mathcal D(\widetilde M).$$
Thus $g_t$ fixes each point in 
$\bar{\mathcal D}(P)$. We deduce that $g_t$ is a trivial one-parameter 
subgroup, which is a contradiction.  $\square$\\\\
{\bf 2.3. Sphere-extension property}

Let $M$ be a $C^\infty$ manifold. Let $\mathcal F_1$ and $\mathcal F_2$ be two continuous 
foliations with $C^1$ leaves on $M$ such that 
$$ T\mathcal F_1\oplus  T\mathcal F_2 =TM.$$
If $\mathcal F_1$ is a foliation by planes, i.e. 
each leaf of $\mathcal F_1$ is $C^1$ diffeomorphic 
to a certain $\mathbb R^n$, then $(\mathcal F_1, \mathcal F_2)$ is said 
to be a {\it plane foliation couple}. The local leaves of $\mathcal F_1$ are natural transverse sections of 
$\mathcal F_2$ and we consider only the holonomy maps of $\mathcal F_2$ with respect to these 
special transverse sections.

For each leaf $F_{1,x}$ of $\mathcal F_1$ we denote by $S_{1, x}$ its one-point compactification 
which is homeomorphic to a standard sphere. The point at infinity of $S_{1, x}$ is denoted by $\infty$.\\\\
{\bf Definition 2.1.} Under the notations above, a plane foliation couple $(\mathcal F_1, \mathcal F_2)$ is 
said to have the {\it sphere-extension} property if for each holonomy 
map $\theta$ of $\mathcal F_2$ sending $x$ to $y$ there 
exists a homeomorphism $\Theta: S_{1, x}\to S_{1, y}$ which coincides 
locally with $\mathcal F_2$-holonomy maps on 
$S_{1, x}\diagdown \{\infty, \Theta^{-1}(\infty)\}$ and extends the germ of $\theta$ at $x$.

If $\bar \phi_t$ is a lifted flow of an $C^\infty$ Anosov flow $\phi_t$, then $\bar\phi_t$ is said to have 
the {\it sphere-extension} property if $(\bar{\mathcal F}^+, \bar{\mathcal F}^{-, 0})$ and 
$(\bar{\mathcal F}^-, \bar{\mathcal F}^{+, 0})$ have both the sphere-extension property. 

$\ $

Recall that $\bar{\mathcal F}^+$ 
and $\bar{\mathcal F}^-$ are both foliations by planes. The corresponding 
notion for Anosov diffeomorphisms is defined similarly.

Denote by $(\widetilde {\mathcal F}_1, \widetilde {\mathcal F}_2)$ the lifted couple on $\widetilde M$ of 
$({\mathcal F}_1, {\mathcal F}_2)$. Then it is easily seen that 
if $(\widetilde {\mathcal F}_1, \widetilde {\mathcal F}_2)$ has the sphere-extension property, 
then $({\mathcal F}_1, {\mathcal F}_2)$ has also this property. So by considering the 
lifted flows and 
drawing pictures, we can easily see that the geodesic flows of closed negatively curved manifolds 
have the sphere-extension property.

It is easily verified that hyperbolic 
infra-nilautomorphisms have the 
sphere-extension property. However by {\bf [Man]} each Anosov diffeomorphism defined on a infra-nilmanifold is 
topologically conjugate to a hyperbolic infra-nilautomorphism. We deduce that 
the suspensions of Anosov diffeomorphisms on infra-nilmanifolds have the sphere-extension property.\\\\
{\bf Lemma 2.5.} {\it Let $\phi_t$ and $\psi_t$ be two $C^1$ orbit equivalent $C^\infty$ Anosov flows. 
Suppose that the strong stable and strong unstable distributions of $\psi_t$ are of dimension at least two. 
If $\psi_t$ has the sphere-extension property, then $\phi_t$ has also this property.}\\\\
{\bf  Proof.} Denote by $\phi$ the $C^1$ orbit equivalence. Define $\hat\phi_t= \phi\circ\phi_t\circ\phi^{-1}$ and denote 
by ${\bar E}^\pm$ the strong distributions of $\psi_t$. Define $\hat{\mathcal F}^\pm= \phi({\mathcal F}^\pm)$ and 
$\hat{\mathcal F}^{\pm,0}= \phi({\mathcal F}^{\pm,0})$. It is clear that $\phi_t$ has the sphere-extension property iff 
$\hat\phi_t$ has this property in the natural sense.

Since $\phi_t$ and $\psi_t$ are $C^1$ orbit equivalent, then $\hat{\mathcal F}^{\pm,0}= \bar {\mathcal F}^{\pm,0}$. 
Take a leaf of ${\mathcal F}^+$ and take a non-periodic point $x$ in this leaf. We can identify 
$\hat W^+_x$ and $\bar W^+_x$ naturally as following. 

For all $y\in \hat W^+_x,$ there exists $ t\in \mathbb{R}$ such that $ \psi_t(y)\in \bar W^+_x.$ 
If $\hat W^{+, 0}_x$ contains no periodic orbit, this number $t$ is unique for each $y$ in $\hat W^+_x.$ 
If $\hat W^{+, 0}_x$ contains a periodic orbit, then it contains exactly one periodic orbit. Denote by $T$ its 
mininal positive period with 
respect to $\psi_t$. So if $ \psi_t(y)\in W^+_x,$ then for all $k\in \mathbb{Z}$, $\psi_{t+ k\cdot T}(y)\in \bar W^+_x.$

Conversely if $ \psi_{t_1}(y)\in\bar W^+_x$ and $ \psi_{t_2}(y)\in \bar 
W^+_x$, then $\psi_{t_2-t_1}\bar W^+_x =\bar W^+_x$. 
Thus $t_2-t_1\in T\cdot\mathbb{Z}.$ So by associating $t+ T\cdot\mathbb{Z}$ to $y$, we get a well-defined 
$C^\infty$ map from $\hat W^+_x$ to $\mathbb{R}/{T\mathbb{Z}}$.

Thus by taking a lift if necessary, there exists a unique $C^1$ map $\theta_x : \hat W^+\to \mathbb{R}$ such that
$$\theta_x(x) =0,\  \psi_{\theta_x(y)}(y)\in \bar W_x^+, \  \forall\  y\in \hat W^+_x.$$
Define a $C^1$ map $\hat\eta_x : \hat W^+\to \bar W_x^+$ such that
$$ \hat\eta_x(y) = \psi_{\theta_x(y)}(y).$$
Then $\hat\eta_x$ is easily seen to be a local $C^1$ diffeomorphism such that $\hat\eta(x)=x.$ 
Similar we get $\bar\eta:\bar W^+_x\to\hat W^+_x$ such that $\bar\eta(x)=x.$ If $\hat W^{+,0}_x$ 
contains no $\psi_t$-periodic orbit, then $\hat{\eta}$ and $\bar{\eta}$ are both $C^1$ diffeomorphisms. 

Suppose that $\bar W^{+,0}_x$ contains a unique $\psi_t$-periodic orbit of period $T$. Denote by $z$ the unique 
intersection point of this periodic orbit with $\bar W^+_x$. For each $k\in\mathbb{Z}$, we define
$$\Lambda_k= \{ y\in \bar W^+_x\diagdown \{z\}\mid \hat{\eta}\circ\bar{\eta}(y)= \psi_{kT}(y)\}.$$
Then $\bar {W}^+_x\diagdown \{z\}$ is the disjoint union of $\{\Lambda_k\}_{k\in\mathbb{Z}}$. 
Each $\Lambda_k$ is closed in $\bar {W}^+_x\diagdown \{z\}$ and $x\in\Lambda_0$. 
Take $y\in\Lambda_0$ and since $\psi_{-T}$ 
is a contracting diffeomorphism of $\bar{W}^+_x$, then a small ball 
containing $y$ intersects with at most finitely many $\Lambda_l$ non-trivially. We deduce that $\Lambda_0$ is open. 

Since $\bar E^+$ is at least two-dimensional, then $\bar {W}^+_x\diagdown \{z\}$ 
is connected . We deduce that $\Lambda_0= \bar {W}^+_x\diagdown \{z\}$, i.e. $\hat\eta\circ\bar\eta=Id$. 
Similarly we have $\bar\eta\circ\hat\eta=Id$. We identify $\hat W^+_x$ and $\bar W^+_x$ under these 
two sliding $C^1$ diffeomorphisms $\bar \eta$ and $\hat\eta$. 
We can identify $\hat W^-_x$ and $\bar W^-_x$ similarly. 

Since these identifications conjugate the holonomy maps and $\psi_t$ has the 
sphere extension property, then $\hat\phi_t$ has also this property. We deduce that $\phi_t$ has the 
sphere-extension property.  $\square$\\\\
{\bf 3. Proofs of Theorems 1.2 and 1.3}\\
{\bf 3.1. Construction of a transverse geometric structure}

Denote by $\phi_t$ a $C^\infty$ topologically transitive 
quasiconformal Anosov flow such that $E^+$ and $E^-$ are of 
dimensions at least three. 
Then by Theorem $2.1$ and Lemma 
$2.2$, $E^{+, 0}$ and $E^{-, 0}$ are both $C^\infty$ and there exist 
$\tau^+$ and $\tau^-$ two continuous 
$\phi_t$-invariant conformal structures on $E^+$ and $E^-$ which are $C^\infty$ along the leaves of 
${\mathcal F}^+$ and ${\mathcal F}^-$.

Denote by $\Phi$ the orbit foliation of $\phi_t$. 
For each transverse section $\Sigma$ of $\Phi$ we get two $C^\infty$ foliations ${\mathcal F}^+_\Sigma$ 
and ${\mathcal F}^-_\Sigma$ on $\Sigma$ by intersecting ${\mathcal F}^{\pm,0}$ with $\Sigma$. Denote their 
tangent distributions by $E^+_\Sigma$ and $E^-_\Sigma$ respectively.

We can identify $E^\pm_\Sigma$ and $E^\pm$ by projecting $E^\pm$ onto $E^\pm_\Sigma$ parallel to the 
direction of the flow. Under this identification, we get two conformal structures $\tau^+_\Sigma$ and 
$\tau^-_\Sigma$ on $E^+_\Sigma$ and $E^-_\Sigma$. Since $\tau^-_\Sigma$ is easily seen to be 
invariant under the $\Phi$-holonomy maps and the ${\mathcal F}^+_\Sigma$-holonomy maps, then $\tau^-_\Sigma$ 
is $C^\infty$ on $\Sigma$. Similarly we can 
see that $\tau^+_\Sigma$ is also $C^\infty$ on $\Sigma$. So we get on each 
transverse section $\Sigma$ a 
$C^\infty$ geometric structure $({\mathcal F}^\pm_\Sigma, \tau^\pm_\Sigma)$ which is invariant under the 
$\Phi$-holonomy maps.

Denote by $c_n$ the canonical conformal structure on the 
$n$-dimensional sphere $S^n$ and by $M_n$ the isometry group of $c_n$. Then 
$M_n$ acts transitively on $S^n$ and is called the 
M$\ddot {\rm o}$bius group. Suppose that $E^+$ is of dimension $n$ and 
$E^-$ is of dimension $m$. Then we can 
construct as following a transverse 
$(M_n\times M_m, 
S^n\times S^m)$-structure on $\Phi$.  

For any $x\in M$ we denote by $\bar S^+_x$ and $\bar S^-_x$ the one-point compactifications of $E_x^+$ and 
$E^-_x$. Then they admit naturally $C^\infty$ conformal structures extending $\sigma_x^+$ and $\sigma_x^-$. Since 
$$(h^+_x)_\ast(\tau^+)=\sigma^+_x\  \  {\rm and}\  \  (h^-_x)_\ast(\tau^-)=\sigma^-_x,$$
then $S^+_x$ and $S^-_x$, i.e. 
the one-point compactifications of 
$W^+_x$ and $W^-_x$ admit also natural conformal structures isometric to those of $\bar S^+_x$ and 
$\bar S^-_x$ under the natural extensions of $h^+_x$ and $h^-_x$, which are 
denoted by $\bar h^+_x$ and $\bar h^-_x$.

By fixing two conformal frames of $E^+_x$ and $E^-_x$ we get two $C^\infty$ conformal isometries 
$\phi_x^+: \bar S^+_x\to S^n$ and $\phi_x^-: \bar S^-_x\to S^m.$

Take a $C^\infty$ small transverse section $\Sigma_x$ containing $x$ and pieces of 
$W^+_x$ and $W^-_x$. Thus for $\delta\ll 1$ we get the local diffeomorphism 
$$\theta_x: W^+_{x, \delta}\times W^-_{x, \delta} \to \Sigma_x$$
$$(y, z)\to W^-_{\Sigma_x, y, 2\delta}\cap W^+_{\Sigma_x, z, 2\delta}.$$
Then we define 
$\phi_x: \Sigma_x\to S^n\times S^m$ such that $\phi_x= (\phi^+_x\times\phi^-_x)\circ (\bar h^+_x\times \bar h^-_x)\circ
\theta_x^{-1}.$ 

Since $\tau^+_{\Sigma}$ and $\tau^-_{\Sigma}$ are invariant under respectively the 
${\mathcal F}^-_\Sigma$-holonomy maps and the ${\mathcal F}^+_\Sigma$-holonomy maps, then by 
its definition, $\phi_x$ is easily seen to be a local 
isometry of $({\mathcal F}^\pm_{\Sigma_x}, \tau^\pm_{\Sigma_x})$ to $(\{S^n\times\ast\}, \{\ast
\times S^m\}, c_n\times c_m)$. 

Let $h$ be any $\Phi$-holonomy map from a germ of $\Sigma_x$ to 
a germ of $\Sigma_y$. Then it is easy to see that $\theta_y^{-1}\circ h\circ\theta_x$ is given by 
weak holonomy maps. We deduce that $\phi_y\circ h
\circ\phi_x^{-1}= \phi\times \psi$, where $\phi$ and $\psi$ are respectively 
local conformal isomertries of $S^n$ and $S^m$. Since $n, m\geq 3$, then by the following classical theorem of 
Liouville, $\phi$ and $\psi$ can be both extended to global conformal isometries of $S^n$ and $S^m$. So 
$\{(\Sigma_x, \phi_x)\}_{x\in M}$ gives a transverse $(M_n\times M_m, S^n\times S^m)$-structure of $\Phi$.\\\\
{\bf Theorem 3.1.} (Liouville) {\it For 
$n\geq 3$, each local conformal isometry of $S^n$ defined on a connected open subset can 
be extended uniquely to a global conformal isometry.}\\\\
{\bf 3.2. Completeness }

Fix a developing map $\mathcal D$ of the transverse $(M_n\times M_m, S^n\times S^m)$-structure 
of $\Phi$ defined in the previous subsection. Denote by $\mathcal H$ the associated holonomy 
representation.\\\\
{\bf Lemma 3.1.} {\it Under the notations above, each leaf of $\widetilde{\mathcal F}^+$ intersects 
each leaf of $\widetilde{\mathcal F}^{-, 0}$ at most once.}\\\\
{\bf Proof.} By the definition of $\mathcal D$, for any 
$x\in\widetilde M$, ${\mathcal D}(\widetilde W^+_{x})$ must 
be contained in a certain subset 
$S^n\times b$. Then for any 
$y\in \widetilde W^+_{x}$, there exists $\gamma_1\times \gamma_2\in M_n\times M_m$ such that 
$$ (\gamma_1\times \gamma_2)\circ\mathcal D= \widetilde\phi_y,$$
where $\widetilde\phi_y$ is defined similarly as above. Denote by 
$\mathcal D_1$ the composition $pr_1\circ \mathcal D$. Then there exists an open neighborhood 
$V_y$ of $y$ in $\widetilde W^+_{x}$ such that 
$$\gamma_1\circ \mathcal D_1\mid_{V_y}=\widetilde\phi^+_y\circ \widetilde h^+_y.$$
Since $\widetilde h^+_x\circ{\widetilde {h^+_y}}^{-1}$ sends 
$\widetilde \sigma^+_y$ to $\widetilde \sigma^+_x$ and the dimension of $E^+$ 
is at least two, then $\widetilde h^+_x\circ{\widetilde {h^+_y}}^{-1}$ is an affine map. Thus there 
exists $\gamma\in M_n$ such that 
$$\gamma\circ \mathcal D_1\mid_{\widetilde W^+_{x}} = \widetilde\phi^+_x\circ \widetilde h^+_x.$$ 
So $\mathcal D$ sends 
$\widetilde W^+_{x}$ diffeomorphically onto a set of the form $(S^n\diagdown a)\times b$. 

For any $y\in S^n\diagdown a$ such that 
$\mathcal D(z)=y$ and $z\in \widetilde W^+_{x}$, $\mathcal D$ sends $\widetilde W^-_{z}$ diffeomorphically onto 
a set of the form $y\times (S^m\diagdown \omega(y))$. So we get a 
well-defined map $\omega: S^n\diagdown a\to S^m$. 

Now suppose that $\widetilde W^+_x$ intersects $\widetilde W^{-, 0}_x$ at a point $x'$ other than $x$. Then 
there exist $y, y'\in S^n$ such that $y\not= y'$ and 
$$\mathcal D(W^-_x)\subseteq y\times S^m,\  \  \mathcal D(W^-_{x'})\subseteq {y'}\times S^m.$$
Denote by $x''$ the intersection of the $\widetilde\phi_t$-orbit of $x'$ with $\widetilde W^-_x$. Then we have 
$\mathcal D(x'')\not= \mathcal D(x')$. However by the definition of $\mathcal D$, 
$\mathcal D(x'') = \mathcal D(x')$, which is a contradiction. We deduce that 
each leaf of $\widetilde{\mathcal F}^+$ intersects 
each leaf of $\widetilde{\mathcal F}^{-, 0}$ at most once.  $\square$

$\ $
 
The following lemma is a direct 
consequence of the previous lemma, which is firstly observed by T. Barbot in {\bf [Ba]}.\\\\
{\bf Lemma 3.2.} {\it Under the notations above, the lifted orbit space $Q_{\widetilde\Phi}$ is Hausdoff.}\\\\
{\bf Proof.} Suppose on the contrary that there exist two different orbits 
$\widetilde \Phi_1$ and $\widetilde\Phi_2$ such that each 
$\widetilde\Phi$-saturated open neighborhood of $\widetilde \Phi_1$ intersects that of $\widetilde\Phi_2$. 
We want to see that these two orbits are contained in the same leaf of $\widetilde{\mathcal F}^{-,0}$.

Suppose that it is not the case. Denote by $F_1$ and $F_2$ the leaves of $\widetilde{\mathcal F}^{-, 0}$ containing 
respectively these two orbits. Then by assumption the $\widetilde{\mathcal F}^+$-saturated sets 
of $F_1$ and $F_2$ intersect non-trivially. We deduce that 
there exists a leaf $\widetilde W^+_x$ intersecting $F_1$ and $F_2$. Denote by $V_1$ and $V_2$ two 
disjoint open subsets of $\widetilde W^+_x$ containing respectively the intersection of $\widetilde W^+_x$ with $F_1$ 
and that of $\widetilde W^+_x$ with $F_2$. Then by assumption the 
$\widetilde{\mathcal F}^{-, 0}$-saturated set of $V_1$ intersects that of $V_2$ non-trivially, which contradicts 
Lemma $3.1$. 

Thus $\widetilde \Phi_1$ and $\widetilde\Phi_2$ are contained in the same leaf of $\widetilde{\mathcal F}^{-,0}$. 
Similar we can prove that they are contained in the same leaf of $\widetilde{\mathcal F}^{+,0}$. Then 
by Lemma $3.1$, we have $\widetilde \Phi_1 =\widetilde\Phi_2$, which is a contradiction.  $\square$

$\ $

For each $x\in \widetilde M$ we construct an open subset $U_x$ of $\widetilde M$ such that 
$U_x$ is the union of the leaves of $\widetilde {\mathcal F}^{-, 0}$ intersecting 
$\widetilde W^{+, 0}_x$. 
Then we can find a sequence $\{x_i\}_{i=1}^\infty\subseteq\widetilde M$ satisfying the following conditions:\\
$(1)$ $\cup_{i\geq 1}U_{x_i}=\widetilde M.$\\
$(2)$ For each $k\geq 1$, $\Omega_k= \cup_{i=1}^k U_{x_i}$ is connected.\\
In the following we denote $U_{x_i}$ by $U_i$. Largely 
inspired by the arguments in {\bf [Gh1]}, we prove the following 
lemma.\\\\
{\bf Lemma 3.3.} {\it For each $k\geq 1$, $\mathcal D\mid_{\Omega_k}: \Omega\to \mathcal D(\Omega_k)$ is a 
$C^\infty$ fiber bundle with fiber $\mathbb R$ over $\mathcal D(\Omega_k)$. In addition $\mathcal D(\Omega_k)$ is 
either the complement in $S^n\times S^m$ of the graph of a continuous map from $S^n$ to $S^m$ or the 
complement of the union of $\{\ast\}\times S^m$ and of the graph of a continuous map from 
$(S^n\diagdown\{\ast\})$ to $S^m$.}\\\\
{\bf Proof.} We prove this lemma by induction. For $k=1$ we have 
$\Omega_1=U_1$. In the proof of 
Lemma $3.1$ we have seen that $\mathcal D$ sends 
$\widetilde W^+_{x}$ diffeomorphically onto a set of the form $(S^n\diagdown a)\times b$. 
For any $y\in S^n\diagdown a$ such that 
$\mathcal D(z)=y$ and $z\in \widetilde W^+_{x_1}$, $\mathcal D$ sends $\widetilde W^-_{z}$ diffeomorphically onto 
a set of the form $y\times (S^m\diagdown \omega(y))$. So we get a 
well-defined map $\omega: S^n\diagdown a\to S^m$. 

Denote by $Gr(\omega)$ the graph of $\omega$. Then the complement of $Gr(\omega)$ in 
$(S^n\diagdown a)\times S^m$ is the open set $\mathcal D(U_1)$. So $\omega$ is continuous. By the definition of 
$U_1$ the inverse images of $\mathcal D\mid_{U_1}$ are all connected. Then by the existence of fine transverse 
sections $\mathcal D\mid_{U_1}$ is seen to be a fiber bundle of fiber $\mathbb R$. So the lemma is true for $k=1$. 

Suppose that the lemma is true for $\Omega_k$. Then $\mathcal D\mid_{\Omega_k}$ is a fiber bundle with 
fiber $\mathbb R$ and $\mathcal D(\Omega_k)$ is the complement in 
$S^n\times S^m$ of the graph of a $C^0$ map $u_k: S^n\to S^m$ or of the union of a vertical 
$a_k\times S^m$ and the graph of a $C^0$ map $u_k: S^n\diagdown a_k\to S^m$. 

In addition by the argument above, we know that $\mathcal D(U_{k+1})$ is the complement of the union of 
$b_{k+1}\times S^m$ and of the graph of a $C^0$ map $v_{k+1}: S^n\diagdown b_{k+1}\to S^m$ and 
$\mathcal D\mid_{U_{k+1}}$ is a fiber bundle with fiber $\mathbb R$.

Since $\mathcal D(U_{k+1})\cap\mathcal D(\Omega_k)$ is the complement in $S^n\times S^m$ of a finite 
union of topological submanifolds of codimension at least two, then 
$\mathcal D(U_{k+1})\cap\mathcal D(\Omega_k)$ is connected and open.

Firstly we want to see that $\mathcal D\mid_{\Omega_{k+1}}$ is a fiber bundle with fiber $\mathbb R$. Take 
$x\in \Omega_k$ and $y\in U_{k+1}$ such that $\mathcal D(x)=\mathcal D(y)$. Since $\Omega_{k+1}$ is 
connected, then $\Omega_k\cap U_{k+1}\not=\emptyset$. So we can take $z\in \Omega_k\cap U_{k+1}$ and a 
$C^0$ curve $\gamma$ in $\mathcal D(U_{k+1})\cap\mathcal D(\Omega_k)$ connecting $\mathcal D(z)$ and 
$\mathcal D(x)$.

Since $\mathcal D\mid_{\Omega_k}$ and $\mathcal D\mid_{U_{k+1}}$ are fiber bundles with 
fiber $\mathbb R$, then we can lift $\gamma$ to 
two $C^0$ curves $\gamma_1\subseteq \Omega_k$ and $\gamma_2\subseteq U_{k+1}$ such that 
$\gamma_1(0)=\gamma_2(0)=z$. Thus $\gamma_1(1)$ 
and $x$ are contained in the same 
$\widetilde\phi_t$-orbit and so are $\gamma_2(1)$ 
and $y$.

Denote by $\Lambda$ the subset of $t\in[0, 1]$ such that 
$\gamma_1(t)$ and $\gamma_2(t)$ are in the same orbit of $\widetilde \phi_t$. By the 
section property, $\Lambda$ is 
easily seen to be open in $[0, 1]$. Suppose 
that $\{t_n\}_{n=1}^\infty\subseteq 
\Lambda$ and $t_n\to t$. If 
$\gamma_1(t) $ and $\gamma_2(t)$ 
are not in the same $\widetilde\phi_t$-orbit, 
then by Lemma $3.2$ there exist disjoint $\widetilde\phi_t$-saturated 
open neighborhoods of $\gamma_1(t) $ and $\gamma_2(t)$. Thus for $n\gg 1$, $\gamma_1(t_n) $ and 
$\gamma_2(t_n)$ are not in the same orbit of $\widetilde\phi_t$, which is a contradiction. We deduce that 
$\Lambda$ is closed. Thus $\Lambda=[0, 1]$. So $x$ and $y$ are contained in the 
same $\widetilde\phi_t$-orbit. 
We deduce that $\mathcal D\mid_{\Omega_{k+1}}$ is a fiber bundle with fiber $\mathbb R$.

Now we want to see the form of $\mathcal D(\Omega_{k+1})$. Suppose at first that 
$\mathcal D(\Omega_{k}) = (Gr(u_k))^c$. Take $p\in S^n\diagdown b_{k+1}$. If $u_k(p)\not= v_{k+1}(p)$, then 
$\mathcal D(\Omega_{k+1})$ contains the vertical $p\times S^m$. Since $p\not =b_{k+1}$, then there exists 
$x\in U_{k+1}$ such that $\mathcal D(\widetilde W^-_x)= p\times(S^m\diagdown v_{k+1}(p))$. So there 
exists $y\in \Omega_k$ such that $ p\times v_{k+1}(p)\in \mathcal D(\widetilde W^-_y)$. 
In particular, $\mathcal D(\widetilde W^-_x)\cap \mathcal D(\widetilde W^-_y)\not= \emptyset$. So there exists 
$t\in \mathbb R$ such that $\widetilde \phi_t(\widetilde W^-_x) =\widetilde W^-_y$. We deduce that 
$\mathcal D(\widetilde W^-_x)= p\times S^m$, which is absurd. 
So in this case, $\mathcal D(\Omega_{k+1}) = (Gr(u_k))^c$.

Suppose that $\mathcal D(\Omega_{k}) = (Gr(u_k)\cup(a_k\times S^m))^c$. 
For each $p\in S^n\diagdown\{a_k, b_{k+1}\}$ 
we get as above that $u_k(p)= v_{k+1}(p)$. 

If $a_k\not= b_{k+1}$, then $u_k$ and $v_{k+1}$ can be extended to 

the same continuous map $\bar u_k$ on $S^n$. In this case $\mathcal D(\Omega_{k+1})= (Gr(\bar u_k))^c.$ 

If $a_k= b_{k+1}$, then we certainly 
have $\mathcal D(\Omega_{k+1})= (Gr(u_k)\cup(a_k\times S^m))^c$.  $\square$

$\ $

We deduce from the previous lemma that $\mathcal D: \widetilde M\to \mathcal D(\widetilde M)$ is a $C^\infty$ 
fiber bundle with fiber $\mathbb R$. So the transverse $(M_n\times M_m, S^n\times S^m)$-structure of $\Phi$ is 
complete. In addition by the proof of the previous lemma we see that 
if $b_k$ is equal to $b_1$ for each $k\geq 1$ then $\mathcal D(\widetilde M)= (Gr(u_1)\cup(a_1\times S^m))^c$. 
If there exists $k>1$ such that $b_k\not =b_1$ then $\mathcal D(\widetilde M)= (Gr(\bar u_1))^c$.

By exchanging the roles of $E^+$ and $E^-$ in the previous lemma, we get the following two cases:\\
$(1)$ $\mathcal D(\widetilde M)= (S^n\diagdown a)\times (S^m\diagdown b).$\\
$(2)$ $\mathcal D(\widetilde M)= (Gr(f))^c$ where $f$ is a homeomorphism of $S^n$ onto $S^m$. In particular 
$n=m$ in this case.

Let us consider firstly Case $(1)$. By changing the developing map we can suppose that $a=b=\infty$. Denote by 
$CO_n$ the isometry group of the canonical conformal structure of ${\mathbb R}^n$. Then by Lemma $2.3$ we get 
a compatible transverse $( CO_n\times CO_m, {\mathbb R}^n\times{\mathbb R}^m)$-structure of $\Phi$. 
In particular, the weak stable and weak unstable foliations admit transverse affine structures. 
So by {\bf [Pl2]} the flow $\phi_t$ 
admits a $C^\infty$ global section $\Sigma$. Since the Poincar\'e map $\phi$ of $\Sigma$ is also 
topologically 
transitive and quasiconformal, then by {\bf [K-Sa]}, $\phi_t$ is $C^\infty$ conjugate to a finite factor of 
a hyperbolic automorphism of a torus. 
We deduce that up to finite covers, $\phi_t$ is $C^\infty$ orbit equivalent to the suspension of a 
hyperbolic automorphism of a torus.

Now we consider Case $(2)$. Denote by $\Gamma$ 
the fundamental group of $M$. Then by Lemma $2.4$ the group $\mathcal H(\Gamma)$ is discrete in 
$M_n\times M_n$. Define $\mathcal H_1=pr_1\circ\mathcal H$ and $\mathcal H_2=pr_2\circ\mathcal H$. Then we 
have 
$$ f\circ\mathcal H_1(\gamma)\circ f^{-1}= \mathcal H_2(\gamma),\  \forall\  \gamma\in \Gamma.$$
We deduce that $\mathcal H_1(\Gamma)$ and $\mathcal H_2(\Gamma)$ are both discrete in $M_n$. Denote 
$\mathcal 
H_1(\Gamma)$ by $\Gamma_1$ and $\mathcal H_2(\Gamma)$ by $\Gamma_2$. Since $\phi_t$ is topologically 
transitive, then $\Phi$ admits at least a simply connected leaf. We deduce that $\mathcal H$ is injective. 
So $\Gamma$ and $\Gamma_1$ and $\Gamma_2$ are all isomorphic.

We can prove that $\Gamma_1$ is uniform in $M_n$ 
as following. Suppose on the contrary that $\Gamma_1$ is not uniform. 
Then $\Gamma_1$ admits a finite index torsion free subgroup $\Gamma'_1$ such that $cd(\Gamma'_1)\leq n$, 
where $cd(\Gamma'_1)$ denotes the cohomological dimension of $\Gamma'_1$. So by passing to a 
finite index subgroup if necessary, we can suppose that $cd(\Gamma)\leq n$. 

Denote by $B\Gamma$ the classifying 
space of $\Gamma$ and by $E\Gamma$ the universal covering space of $B\Gamma$. Then we have
$$ B\Gamma\cong \Gamma_1\diagdown\mathbb H^{n+1},\  E\Gamma\cong \mathbb H^{n+1},$$
where ${\mathbb H}^{n+1}$ denotes the simply connected hyperbolic space of dimension $n+1$. Denote by 
$E\Gamma\times_\Gamma\widetilde M$ the quotient manifold of $E\Gamma\times\widetilde M$ under the diagonal 
action of $\Gamma$. Then we have the following fibre bundle with fiber $\widetilde M$ 
$$\pi_1: E\Gamma\times_\Gamma\widetilde M\to B\Gamma,$$
$$\Gamma((a, x))\to \Gamma(a).$$
By using the cohomology Leray-Serre spectral sequence to this fibre bundle (see {\bf [Mc]}), we get that 
$$E^{p, q}_2= H^p(\Gamma, H^q(\widetilde M))$$ 
converges to $H^{p+q}(E\Gamma\times_\Gamma\widetilde M)$.
Since $\widetilde M$ is a fibre bundle with fiber $\mathbb R$ and base 
$(S^n\times S^n)\diagdown (Gr(f))$, then $\widetilde M$ 
is homotopically equivalent to the sphere $S^n$. Since we 
have in addition $cd(\Gamma)\leq n$, then we deduce from the spectral 
sequence above that $H^{2n+1}(E\Gamma\times_\Gamma\widetilde M)$ is trivial. 

However by projecting onto the second factor $E\Gamma\times_\Gamma\widetilde M$ is easily seen to be also 
a fibre bundle over $M$ and with contractible fiber $E\Gamma$. So $E\Gamma\times_\Gamma\widetilde M$ is 
homotopically equivalent to $M$. We deduce that $H^{2n+1}(M)$ is trivial, which is absurd. 
So $\Gamma_1$ is uniform in $M_n$. Similarly $\Gamma_2$ is also uniform in $M_n$. 

Since $f$ conjugates $\Gamma_1$ to $\Gamma_2$, then by Mostow's rigidity theorem 
(see {\bf [Mo]}) $f$ is contained in $M_n$. 
So by replacing $\mathcal D$ by $(Id\times f^{-1})\circ\mathcal D$, we can suppose that $f=Id$ and 
$$\mathcal D(\widetilde M)= (S^n\times S^n)\diagdown\Delta,$$
where $\Delta$ denotes the diagonal of $S^n\times S^n$. In addition, 
we have $\mathcal H_1= \mathcal H_2$. So by Lemma $6.3.1$, $\Phi$ admits a 
compatible transverse $(M_n, (S^n\times S^n)\diagdown\Delta)$-structure with respect to the 
diagonal action of $M_n$ on $(S^n\times S^n)\diagdown\Delta.$ 

Lift $\phi_t$ to a finite cover to eliminate the torsion of $\Gamma$ and define 
$V=\mathcal H (\Gamma)\diagdown\mathbb{H}^{n+1}$. 
Then $V$ is a closed hyperbolic manifold. In 
addition, the $\Gamma$-action on $Q_{\widetilde \Phi}$ is $C^\infty$ conjugate to the 
$\mathcal H(\Gamma)$-action on the leaf space of the lifted geodesic flow of $V$ under 
$\mathcal D$ and $\mathcal H$. Since the holonomy of each periodic orbit of 
$\phi_t$ is non-trivial, then the holonomy covering of each leaf of $\Phi$ is contractible. Denote 
by $\psi_t$ the geodesic flow of $V$. So by Proposition $2.1$ 
there exists a $C^\infty$ homotopy equivalence $h$ conjugating the leaf space of $\phi_t$ 
with that of $\psi_t$. 
However $h$ is not in general a $C^\infty$ diffeomorphism. In order to get 
a $C^\infty$ orbit conjugacy between $\phi_t$ and $\psi_t$, we use a 
classical diffusion argument discovered by \'E. Ghys. Let 
us recall briefly this argument (see {\bf [Gh2]} and {\bf [Ba]} for details): 

There exists a $C^\infty$ function $u: \mathbb R\times M\to\mathbb R$ such that 
$$h(\phi_t(x))= \psi_{u(t, x)}(h(x)),\  \forall \  t\in\mathbb R,\  \forall\  x\in M.$$
Define for $T\gg 1$, $u_T(x)=\frac{1}{T}\int_0^T u(s, x) ds$ and 
$h_T: M\to T^1V$ such that $h_T(x)=\psi_{u_T(x)}(h(x))$. If $T\gg 1$, then we can see 
that $h_T$ satisfies the same conditions as $h$ and is a $C^\infty$ diffeomorphism. 

So up to finite covers, $\phi_t$ is $C^\infty$ orbit equivalent to the geodesic flow of a closed hyperbolic manifold, 
which finishes the proof of the first part of Theorem $1.2$.\\\\
{\bf 3.3. Smoothness blowing up}

In this subsection we prove the second part of Theorem $1.2$. Suppose that $\phi_t$ satisifes the conditions of 
Theorem $1.2$ such that $E^+\oplus E^-$ is in addition $C^1$. 
Then because of the first part of Theorem $1.2$, $\phi_t$ 
is seen to be volume-preserving. So in order to prove the second part of Theorem $1.2$, 
we need only prove the $C^\infty$ 
smoothness of $E^+\oplus E^-$ and then use the following classification result established in {\bf [Fa1]} :\\\\ 
{\bf Theorem 3.2.} ({\bf [Fa1]}, Theorem 1) {\it Let $\phi_t$ be a $C^\infty$ volume-preserving 
uniformly quasiconformal Anosov flow on a closed manifold $M$. If 
$E^+\oplus E^-$ is $C^\infty$ and the dimensions of $E^+$ and $E^-$ are at least 2, then 
up to a constant change of time scale 
and finite covers, $\phi_t$ is $C^\infty$ flow equivalent either to the suspension of 
a hyperbolic automorphism of a torus, or to a canonical perturbation of the geodesic flow 
of a hyperbolic manifold.}\\\\
%which classifies the volume-preserving 
%quasiconformal Anosov flows whose strong stable and unstable distributions 
%are of dimension at least $2$ and $E^+\oplus E^-$ is $C^\infty$.\\\\
{\bf Lemma 3.4.} {\it Under the notations above, $E^+\oplus E^-$ is $C^\infty$.}\\\\
{\bf Proof.} Suppose at first that 
$\phi_t$ is $C^\infty$ orbit equivalent to the geodesic flow $\psi_t$ of a hyperbolic manifold (up to finite covers). 
Denote by $\lambda$ the canonical $1$-form of $\phi_t$ and by $X$ the generator of $\psi_t$. 
Up to $C^\infty$ flow conjugacy we suppose that $\phi_t$ is generated by $fX$. 

Since $E^+\oplus E^-$ is supposed to be $C^1$, then $\lambda$ is $C^1$ and 
$\lambda(X)=\frac{1}{f}$ is $C^\infty$. It is easily seen that $d\lambda$ is $\phi_t$-invariant. Then by 
the Anosov property of $\phi_t$, we get $i_{fX}d\lambda=0$. Thus $i_Xd\lambda=0$. We deduce that 
$d\lambda$ is $\psi_t$-invariant. Denote by $\lambda'$ the canonical $1$-form of $\psi_t$. Then by 
{\bf [Ham1]} there exists $a\in\mathbb R$ such that 
$d\lambda= a\cdot d\lambda'$. 

Define $\beta= \lambda- a\cdot\lambda'$. Then $\beta$ is a $C^1$ $1$-form such that 
$d\beta=0$. In addition $\beta(X)= \lambda(X)-a$ is $C^\infty$. \\\\
{\bf Sublemma.} {\it Let $\phi_t$ be a $C^\infty$ volume-preserving Anosov flow 
on $M$. If $\alpha$ is a $C^1$ $1$-form on $M$ such that $d\alpha=0$ and 
$\alpha(X)$ is $C^\infty$, then $\alpha$ is $C^\infty$. }\\\\
{\bf Proof.} Since $d\alpha=0$ and 
the Stokes formula is valid for $C^1$ forms (even for Lipchitz forms), then 
there exists a $C^\infty$ $1$-form $\beta$ giving the same element of 
$(H_1(M, \mathbb R))^\ast$ as that of $\alpha$. So by integrating 
$(\alpha-\beta)$ along curves, we get a well-defined $C^2$ function $f$ on $M$. 
Thus for any $x\in M$ and any $t\in\mathbb R$ we have 
$$f(\phi_t(x))- f(x) =\int_0^t(\alpha(X)-\beta(X))\circ\phi_s(x) ds.$$
Since $\alpha(X)$ is supposed to be $C^\infty$, then by {\bf [LMM]}, $f$ is seen to be 
$C^\infty$. However by the definition of $f$, we have $\alpha-\beta= df$. Thus $\alpha$ is $C^\infty$.  $\square$

$\ $

We deduce from this sublemma that $\beta$ is $C^\infty$. Thus $\lambda$ is $C^\infty$. So 
$E^+\oplus E^-$ is also $C^\infty$. 

If $\phi_t$ is $C^\infty$ orbit equivalent to the suspension $\psi_t$ of a hyperbolic 
automorphism of a torus (up to finite covers), then by similar arguments as above, we 
can see that $d\lambda$ is 
$\psi_t$-invariant.

Take a leaf $\Sigma$ of the foliation of the sum of the strong stable and the strong unstable distributions 
of $\psi_t$ and denote by $\psi$ its Poincar\'e map. Then 
$\lambda\mid_\Sigma$ is $C^1$ and $d(\lambda\mid_\Sigma)$ is $\psi$-invariant. Thus by the same 
arguments as in {\bf [Fa2]} we get 
$d(\lambda\mid_\Sigma)=0$. We deduce that $d\lambda=0$. Since in addition $\lambda(X)=\frac{1}{f}$ is $C^\infty$, 
then by the previous sublemma $\lambda$ is $C^\infty$. Thus $E^+\oplus E^-$ is also $C^\infty$.  $\square$\\\\
{\bf Proof of Theorem 1.3.} Suppose that $\phi_t$ satisfies the conditions of 
Theorem $1.3$. Similar to the previous section, we can construct 
a $C^\infty$ geometric structure $(\mathcal F^\pm_\Sigma, \tau^\pm_\Sigma)$ on each transverse section $\Sigma$ of 
$\Phi$. Similarly we can construct a family of transverse charts $\{(\Sigma_x, \phi_x)\}_{x\in M}$. Then because of the 
sphere-extension property, the chart changes of these charts are easily seen to be 
given by the restrictions of the elements of $M_2\times M_m$ 
with respect to the natural action of $M_2\times M_m$ on $S^2\times S^m$. So in this 
way we get a transverse $(M_2\times M_m, S^2\times S^m)$-structure on $\Phi$. Then as in the previous subsection 
the proof splitts into Case $(1)$ and Case $(2)$. Each of them is understood in the same manner 
as in the previous subsection.\\\\
{\bf 4. Applications to the geodesic flows of hyperbolic manifolds}

Now let us begin to prove Theorem $1.6$.\\\\
{\bf Lemma 4.1.} {\it Let $\phi_t$ and $\psi_t$ be two $C^\infty$ Anosov flows which are $C^1$ orbit equivalent. 
If $\psi_t$ is volume-preserving, then so is $\phi_t$.}\\\\
{\bf  Proof.} By conjugating $\phi_t$ by the $C^1$ orbit conjugacy, we can suppose that $\phi_t$ is a $C^1$ flow 
and a time change of $\psi_t$. Denote by $\nu$ the $\psi_t$-invariant volume form and 
by $X$ the generator of $\psi_t$. Then by taking $i_X\nu$ we get a family of $\Psi$-holonomy invariant volume 
forms 
on the transverse sections 
of $\Psi$. This family of transversal volume forms is also $\Phi$-holonomy invariant. Denote 
by $dt_\phi$ the normalized foliated measure along the leaves of $\Phi$ such that 
$dt_\phi(Y)\equiv 1$, where $Y$ denotes the 
generator of $\phi_t$. In each flow box of $\phi_t$ we 
take the product 
measure $\nu_\Sigma\otimes dt_\phi$. Then it is easily seen that in the intersection of two flow boxes the 
two measures coincide. Then we can extend this family of local measures to a measure 
$\mu$ on $M$ which is in the Lebesgue class and easily seen to be $\phi_t$-invariant.  $\square$\\\\
{\bf  Proof of Theorem 1.6.} Since $\psi_t$ is conformal, then by Lemma $1.1$, $\phi_t$ is quasiconformal. 
In addition by lemmas $2.5$ and 
$4.1$, $\phi_t$ verifies the conditions of Theorem $1.2$ or Theorem $1.3$. 
So up to finite covers, $\phi_t$ is $C^\infty$ orbit equivalent either to a 
suspension or to the geodesic flow of a hyperbolic manifold $\hat\psi_t$. Since $\psi_t$ is contact, then it admits 
no $C^1$ global section. So up to finite covers, $\phi_t$ is $C^\infty$ orbit equivalent to $\hat\psi_t$. 

However in the proofs of Theorems $1.2$ and $1.3$, we passed to a finite cover only in order to eliminate the 
torsion in the fundamental group of $M$. But in the current 
case, the fundamental group has no torsion by the classical Cartan theorem. 
So $\phi_t$ is $C^\infty$ orbit equivalent to $\hat\psi_t$. Then by Mostow's rigidity theorem (see {\bf [Mo]} and 
{\bf [M]}), 
$\hat\psi_t$ is $C^\infty$ flow equivalent to $\psi_t$. We deduce that $\phi_t$ is $C^\infty$ orbit equivalent to 
$\psi_t$.  $\square$\\\\
{\bf Proof of Corollary 1.1.} Let us prove firstly $(1)$. Suppose 
that the geodesic flow of $M$ is H$\ddot{\rm o}$lder-continuously orbit equivalent to that of a hyperbolic 
manifold $N$. Since $n\geq 3$, then the fundamental group of $M$ is isomorphic to that of $N$. Since $\pi_1(M)$ 
with its word metric is quasi-isometric to $\widetilde M$ and $\pi_1(N)$ is quasi-isometric to 
$\mathbb H^n$, then we deduce that $\widetilde M$ is quasi-isometric to $\mathbb H^n$. 

Conversely, if $\widetilde M$ is quasi-isometric to $\mathbb H^n$, then $\pi_1(M)$ is also quasi-isometric to 
$\mathbb H^n$. Thus by {\bf [Su]} and {\bf [Tu]}, there exists a uniform lattice $\Gamma$ in the isometric group of 
$\mathbb H^n$ and a surjective group homomorphism $\rho: \pi_1(M)\to \Gamma$ such that the kernal of 
$\rho$ is finite. However by a classical result of \'E. Cartan, $\pi_1(M)$ is without torsion. We deduce that 
$\pi_1(M)$ is isomorphic to $\Gamma$. In particular, $\Gamma$ is also without torsion. So 
$N= \Gamma\diagdown \mathbb H^n$ is a closed hyperbolic manifold. 

Denote respectively by $\phi_t$ and $\psi_t$ the geodesic flows of 
$M$ and $N$. Since $\pi_1(M)\cong\pi_1(N)$, then by {\bf [Gr]}, $\phi_t$ is $C^0$ orbit equivalent to 
$\psi_t$. Since each continuous 
orbit conjugacy between 
Anosov flows can be $C^0$ approximated by H$\ddot{\rm o}$lder-continuous 
orbit conjugacies (see {\bf [HK]}), then $(1)$ is true. 

Now let us prove $(2)$. We need only prove the necessarity. 
Suppose that the geodesic flow $\phi_t$ of $M$ is $C^1$ 
orbit equivalent to the geodesic flow $\psi_t$ of a closed 
hyperbolic manifold $N$. Since $\psi_t$ is 
conformal, then by Lemma $1.1$, $\phi_t$ is 
quasiconformal. Thus by Corollary $1.1$, it is $C^\infty$ orbit equivalent to the geodesic flow of $N$. 
Since $C^\infty$ orbit conjugacy preserves weak stable and weak unstable distributions, then 
$\phi_t$ is Anosov-smooth. So by {\bf [BFL]}, it is $C^\infty$ flow equivalent to the geodesic flow of a 
Riemannian manifold of 
constant negative curvature. Then by {\bf [BCG]}, $M$ has constant negative curvature. $\square$

$\ $

Before the proof of Corollary $1.2$, let us recall firstly some notions. A $C^\infty$ flow $\phi_t$ defined on a closed 
$n$-dimensional manifold is said to be {\it quasi-Anosov} if there exists a continuous $(n-1)$-dimensional 
distribution $\nu$ transversal to the flow, such that for any non-zero vector $v$ in $\nu$, the set 
$\{\parallel D\phi_t(v)\parallel,\  t\in \mathbb R\}$ is unbounded with 
respect to a certain (then all) Riemannian metric. In {\bf [Ma]}, R. Man$\tilde{\rm e}$ proved the 
following important result.\\\\
{\bf Theorem 4.1.} (R. Man$\tilde{\rm e}$) {\it If $\phi_t$ is quasi-Anosov and volume-preserving, then 
$\phi_t$ is Anosov.}\\\\
{\bf Proof of Corollary 1.2.} Denote $h:M\to N$ a $C^1$ conjuguacy sending the leaves of $\Phi$ onto those 
of $\Psi$. Since $\Phi$ is the orbit foliation of the geodesic flow of a hyperbolic manifold, then $\Psi$ is 
orientable. Thus we can find a $C^\infty$ non-where vanishing vector field $Y$ tangent to $\Psi$, whose flow is 
$C^1$ orbit equivalent to $\phi_t$ under $h$. We denote by $\psi_t$ the flow of $Y$.

There exists a $C^1$ map $\alpha: M\times \mathbb R\to \mathbb R$ such that 
$\psi_t\circ h= h\circ \phi_{\alpha(t,\cdot)}$. Since 
$h$ is $C^1$, then there exists $A>0$ such that for any $u\in TM$, we have
$$ \frac {1}{A}\parallel u\parallel \leq \parallel Dh(u)\parallel\leq A\parallel u\parallel.$$
Denote by $E^+$ and $E^-$ the strong unstable and stable 
distributions of $\phi_t$. Then for any $x\in M$ and any $(u^++ u^-)\in (E^+\oplus E^-)_x $, 
$$\parallel D\psi_t(Dh(u^++u^-))\parallel = \parallel Dh(D(\phi_{\alpha(t, \cdot)})(u^++u^-))\parallel$$
$$\geq \frac{1}{A} \parallel D(\phi_{\alpha(t, x)})(u^+)+ D(\phi_{\alpha(t, x)})(u^-)+ b(t)\cdot 
X_{\phi_{\alpha(t,x)}}(x)\parallel,$$
where $X$ denotes the generator of $\phi$. Since $\phi_t$ is Anosov, then it is easy to see that 
for any $v\in Dh(E^+\oplus E^-)$, $\{\parallel D\psi_t(v)\parallel,\  t\in \mathbb R\}$ is unbounded. 
Thus $\psi_t$ is quasi-Anosov. In addition, we know by Lemma $4.1$ that $\psi_t$ is volume-preserving. Thus 
we deduce from Theorem $4.1$ that $\psi_t$ is a $C^\infty$ Anosov flow, which is $C^1$ orbit equivalent to 
$\phi_t$. Then we deduce from Theorem $1.6$ that $\psi_t$ is $C^\infty$ orbit equivalent to $\phi_t$. Thus 
their orbit foliations $\Psi$ and $\Phi$ are $C^\infty$ conjuguate.\\\\  
%{\bf Remark 5.1.} We can certainly deduce Case $(2)$ of Corollary $1.2$ from the main results in {\bf [Yu]} or 
%{\bf [Sa]}. But it seems more natural for us to deduce it from Corollary $1.1$.\\\\
{\bf Acknowledgements.} The author would like to thank his thesis advisers, P. Foulon and P. Pansu, for the 
discussions and help. He would like also to thank \'E. Ghys and 
B. Hasselblatt for encouragements and T. Barbot and 
J. Lannes for answering questions.  

$\ $

{\bf References}

$\ $

{\small
{\bf [An]} V. D. Anosov, Geodesic flows on closed Riemannian manifolds with negative curvature, 
{\it Proc. Inst. Steklov} 90 (1967) 1-235.

{\bf [Ba]} T. Barbot, Caract\'erisation des flots d'Anosov en dimension 3 par leurs feuilletages faibles, 
{\it Ergod. Th. and Dynam. Sys.} 15 (1995) 247-270.

{\bf [BCG]} G. Besson, G. Courtois and S. Gallot, Entropies et rigidit\'es des espaces localement sym\'etriques 
de courbure strictement n\'egative, {\it GAFA} 5 (1995) 731-799.

{\bf [BFL]} Y. Benoist, P. Foulon and F. Labourie, 
Flots d'Anosov \`a distributions stable et instable diff\'erentiables, 
{\it J. Amer. Math. Soc.} 5 (1992) 33-74.

{\bf [Fa1]} Y. Fang, Smooth rigidity of uniformly quasiconformal Anosov flows, to appear in 
{\it Ergod. Th. and Dynam. Sys.}

{\bf [Fa2]} Y. Fang, A remark about hyperbolic infranilautomorphisms, 
{\it C. R. Acad. Sci. Paris, Ser. I 336 No.9} (2003) 769-772.

{\bf [Gh1]} \'E. Ghys, Holomorphic Anosov flows, {\it Invent. math.} 119 (1995) 585-614.

{\bf [Gh2]} \'E. Ghys, D\'eformation des flots d'Anosov et de groupes fuchsiens, {\it Ann. Inst. Fourier} 42 (1992) 
209-247.

{\bf [Go]} C. Godbillon, Feuilletages, {\it Progress in Mathematics} 98 (1991).

{\bf [Gr]} M. Gromov, Three remarks on geodesic dynamics and fundamental groups, {\it 
Enseign. Math. (2)} 46 (2000) 391-402. 

{\bf [Hae]} A. Haefliger, Groupoides d'holonomie et classifiants, {\it Ast\'erisque} 116 (1984) 70-97.

{\bf [Ham1]} U. Hamenst$\ddot{\rm a}$dt, invariant two-forms for geodesic flows, {\it Math. Ann.} 301 (1995) 
{\it No.4} 677-698.

{\bf [Ham2]} U. Hamenstadt, Cocycles, symplectic structures and intersection, {\it Geom. Funct. Anal.} 9 (1999) 
no.1, 90-140.

{\bf [HK]} B. Hasselblatt and A. Katok, Introduction to the modern theory of dynamical 
systems, {\it Encyclopedia of Mathematics and its Applications, vol 54.} 1995. 

{\bf [K-Sa]} B. Kalinin and V. Sadovskaya, On local and global rigidity of quasiconformal Anosov 
diffeomorphisms, {\it J. Inst. Math. Jussieu No.4} 2 (2003) 567-582.

{\bf [L1]} R. de La llave, Rigidity of higher-dimensional conformal Anosov systems, 
{\it Ergod. Th. and Dynam. Sys.} 22 (2002) no.6, 1845-1870.

{\bf [L2]} R. de La llave, Further rigidity properties of conformal Anosov systems, to appear in 
{\it Ergodic Theory and Dynamical Systems}. 

{\bf [LM]} R. de La llave and R. Moriy\'on, Invariants for smooth conjugacy of hyperbolic dynamical systems, 
IV, {\it Comm. Math. Phys.} 116 (1988) no.2, 185-192.

{\bf [LMM]} R. de la llave, J. Marco and R. Moriyon, Canonical perturbation theory of 
Anosov systems and regularity results for Livsic cohomology equation, {\it Ann. Math} 
123(3) (1986) 537-612.

{\bf [M]} G. A. Margulis, The isometry of closed manifolds of constant negative curvature 
with the same fundamental group, {\it Soviet Math. Dokl.} 11 (1970) 722-723.

{\bf [Ma]} R. Man$\tilde{\rm e}$, Quasi-Anosov diffeomorphisms and hyperbolic manifolds, {\it 
Trans. Amer. Math. Soc.} 229 (1977) 351-369.

{\bf [Man]} A. Manning, There are no new Anosov diffeomorpisms on tori, {\it Amer. J. Math.} 96 (1974) 422-427.  

{\bf [Mc]} J. McCleary, User's Guide to spectral sequences, {\it Mathematics Lectures Series No. 12} (1985).

{\bf [Mo]} G. D. Mostow, Quasi-conformal mappings in $n$-space and the rigidity of hyperbolic space forms, {\it 
Publ. Math. IHES} 34 (1968) 53-104.

{\bf [Pl1]} J. F. Plante, Anosov flows, {\it Amer. J. Math.} 94 (1972) 729-754.

{\bf [Pl2]} J. Plante, Anosov flows, transversely affine foliations and a conjecture of Verjovsky, {\it J. London 
Math. Soc. (2)} 23 (1981) 359-362. 

{\bf [Sa]} V. Sadovskaya, On uniformly quasiconformal Anosov systems, to appear in {\it Math. Res. Lett}.

{\bf [Su]} D. Sullivan, On the ergodic theory at infinity of an arbitrary discrete group of hyperbolic motions, in 
``Riemann surfaces and related topics'', {\it Annals of Math. Studies} 97 (1981) 465-497.

{\bf [Tu]} P. Tukia, On quasiconformal groups, {\it J. Analyse Math.} 46 (1986) 318-346. 

{\bf [Yu]} C. Yue, Quasiconformality in the geodesic flow of negatively curved manifolds. 
{\it GAFA } 6(4) (1996) 740-750.}

\end{document}